\documentclass[11pt,reqno]{amsart}

\usepackage{amsfonts,amscd}
\usepackage{amssymb}
\usepackage{url}
\usepackage[english]{babel}

\theoremstyle{plain}
\newtheorem{theorem}                {Theorem}      [section]
\newtheorem{proposition}  [theorem]  {Proposition}

\theoremstyle{definition}
\newtheorem{example}      [theorem]  {Example}
\newtheorem{remark}       [theorem]  {Remark}

\numberwithin{equation}{section}

\def \R{{\mathbb R}}

\def \C{{\mathbb C}}
\def \H{{\mathbb H}}

\usepackage{color}

\DeclareMathOperator{\grad}{grad}
\DeclareMathOperator{\trace}{trace}

\DeclareMathOperator{\ricci}{Ricci}

\numberwithin{equation}{section}

\begin{document}

\title[]{On  cohomogeneity one biharmonic hypersurfaces into the Euclidean space}

\author{S.~Montaldo}
\address{Universit\`a degli Studi di Cagliari\\
Dipartimento di Matematica e Informatica\\
Via Ospedale 72\\
09124 Cagliari, Italia}
\email{montaldo@unica.it}

\author{C.~Oniciuc}
\address{Faculty of Mathematics\\ ``Al.I. Cuza'' University of Iasi\\
Bd. Carol I no. 11 \\
700506 Iasi, ROMANIA}
\email{oniciucc@uaic.ro}

\author{A.~Ratto}
\address{Universit\`a degli Studi di Cagliari\\
Dipartimento di Matematica e Informatica\\
Viale Merello 93\\
09123 Cagliari, Italia}
\email{rattoa@unica.it}

\begin{abstract}
The aim of this paper is to prove that there exists no cohomogeneity one $G-$invariant proper biharmonic hypersurface into the Euclidean space $\R^n$, where $G$ denotes a tranformation group which acts on $\R^n$ by isometries, with codimension two principal orbits.
This result may be considered in the context of the Chen conjecture, since this family of hypersurfaces includes examples with up to seven distinct principal curvatures. The paper uses the methods of equivariant differential geometry. In particular, the technique of proof  provides a unified treatment for all these $G-$actions.
\end{abstract}

\subjclass[2000]{53A07, 53C42, 58E20}

\keywords{Biharmonic maps, biharmonic immersions, transformation groups, equivariant differential geometry, cohomogeneity one hypersurfaces}

\thanks{Work supported by: PRIN 2010/11 -- Variet\`a reali e complesse: geometria, topologia e analisi armonica N. 2010NNBZ78 003 -- Italy; GNSAGA -- INdAM, Italy;
Romanian National Authority for Scientific Research, CNCS -- UEFISCDI, project number PN-II-RU-TE-2011-3-0108}

\maketitle

\section{Introduction}\label{intro}
According to B.-Y. Chen \cite{Chen} an immersion $\varphi:M^m\hookrightarrow\R^n$ is called {\it biharmonic} if
\begin{equation}\label{def-chen}
\Delta H =(\Delta H_1,\ldots, \Delta H_n)=0 \,\, ,
\end{equation}
where $H=(H_1,\ldots,H_n)$ is the mean curvature vector field and $\Delta$ denotes the Beltrami-Laplace operator on $M$ (our sign convention is such that $\Delta h=-h''$ when $h$ is a function of one real variable).
It follows from Beltrami's equation
$$
m\, H=-(\Delta\, \varphi_1,\ldots,\Delta\, \varphi_n)
$$
that the biharmonicity condition is equivalent to
$$
\Delta^{2}\, \varphi=(\Delta^{2}\, \varphi_1,\ldots,\Delta^{2}\, \varphi_n)=0\,\,,
$$
which justifies the previous definition of biharmonic immersions.

The study of biharmonic immersions in $\R^n$ can be set in a more general variational, Riemannian geometric context. More precisely, we recall that
a smooth map $\varphi:(M,g)\to(N,h)$ is a {\it harmonic map}  if it is a critical point of the {\em energy} functional
\begin{equation}\label{energia}
E(\varphi)=\frac{1}{2}\int_{M}\,|d\varphi|^2\,dv_g \,\, ,
\end{equation}
whose Euler-Lagrange equation is $\tau(\varphi)={\trace} \, \nabla d \varphi =0$.
A natural generalization of harmonic maps are the so-called {\it biharmonic maps}: these maps are the critical points of the bienergy functional (as suggested by Eells--Lemaire \cite{EL83})
\begin{equation}\label{bienergia}
    E_2(\varphi)=\frac{1}{2}\int_{M}\,|\tau (\varphi)|^2\,dv_g \,\, .
\end{equation}
In \cite{Jiang} G.~Jiang showed that the Euler-Lagrange equation associated to $E_2(\varphi)$ is given by $\tau_2(\varphi) =0$, where the {\it bitension field}
$\tau_2(\varphi)$ is
\begin{equation}\label{bitensionfield}
 \tau_2(\varphi) = - \Delta \tau(\varphi)- \trace R^N(d \varphi, \tau(\varphi)) d \varphi\,\,.
\end{equation}
 An immersed submanifold into a Riemannian manifold $(N,h)$ is called a {\it biharmonic submanifold} if the immersion is a biharmonic map.  In particular, minimal immersions are trivially biharmonic, so that we call {\it proper} biharmonic any biharmonic immersion which is not minimal. We observe that when $N=\R^n$ the curvature term in \eqref{bitensionfield} vanishes and $\tau_2(\varphi) =0$ is equivalent to \eqref{def-chen}. Thus the definition of biharmonic immersed submanifolds extends the original definition of Chen. \\

In the case of hypersurfaces, biharmonicity can be expressed by means of the following general result (see \cite{BMO13,C84,LM08,O10,O02}):

\begin{theorem}\label{th: bih subm N}
Let  $\varphi:M^{n-1}\to N^{n}$ be an isometric immersion with mean curvature vector field $H=(f/(n-1))\,\eta$. Then $\varphi$ is biharmonic if and only if the normal and the tangent components of $\tau_2(\varphi)$  vanish,  i.e.,
\begin{subequations}
\begin{equation}\label{eq: caract_bih_normal}
\Delta {f}+f |A|^2- f\ricci^N(\eta,\eta)=0
\end{equation}
and
\begin{eqnarray}\label{eq: caract_bih_tangent}
2 A(\grad f)+  f \grad f-2 f \ricci^N(\eta)^{\top}=0
\end{eqnarray}
\end{subequations}
respectively, where $A$ is the shape operator and $\ricci^N(\eta)^{\top}$ is the tangent  component of the Ricci tensor field of $N$ in the direction of the unit normal vector field $\eta$ of $M$ in $N$.
\end{theorem}

We shall work in the framework of equivariant differential geometry, so let $N$ be a Riemannian manifold and $I(N)$ its full isometry group. It is well-known (see \cite{MST}) that $I(N)$ is a Lie group which acts differentiably on $N$. A Lie subgroup $G$ of $I(N)$ is called an {\it isometry group} of $N$ and, following \cite{HL}, we recall that its {\it cohomogeneity} is defined as the codimension in $N$ of the maximal dimensional orbits, also called the {\it principal} orbits (of course, all the orbits are homogeneous spaces, since they are of the type $G \, / \,H$, where $H$ is the stabilizer). The \emph{cohomogeneity of a G-invariant submanifold} $M$ of $N$ is defined as the dimension of $M$ minus the dimension of the principal orbits.

In this paper we shall be interested in the case that $N$ is the Euclidean space $\R^n$ and the cohomogeneity of $G$ is two, so that \emph{$G-$invariant hypersurfaces are cohomogeneity one submanifolds}. This type of isometry groups of $\R^n$ have been fully classified in \cite{HL}. In particular, Hsiang and Lawson, developping the work of various famous authors, including Cartan and Weil, divided the cohomogeneity two isometry groups $G$ acting on $\R^n$ into five types according to the geometric shape of their orbit space $Q=\R^n \, / \,G$, which is a linear cone in $\R^2$ of angle $\pi \slash d$, $d=1,\,2,\,3,\,4,\,6$ respectively. In this context, a $G-$invariant hypersurface in $\R^n$ can be completely described by means of its profile curve $\gamma$ into the orbit space $Q$. In particular, it turns out that a $G-$invariant hypersurface is a biharmonic submanifold if and only if the curve $\gamma$ satisfies a certain system of ordinary differential equations (see Proposition~\ref{equazioniridottedibiarmonicity} in Section~\ref{equiv-diff-geom-section} below for details). This approach, which uses the symmetries deriving by the group action to reduce a PDE's problem to an ODE's system, has been very fruitful in various context (construction of harmonic maps, counterexamples for Bernstein's type problems for minimal and CMC immersions (see, for example, \cite{BDG,ER, Hsiang}). In general, reduction to an ODE has been a valuable tool because it has helped to produce new solutions. By contrast, in our case what we obtain is a nonexistence result in the direction of the still open Chen conjecture (see \cite{Chen04,Chen} and \cite[chapter~7]{Chen-book}): {\it biharmonic submanifolds into $\R^n$ are minimal}.

Chen's conjecture is still open even for biharmonic hypersurfaces in $\R^n$ though, by a result of Dimitric (see \cite{dimitric}), we know that any biharmonic hypersurface in $\R^n$ with at most two distinct principal curvatures is minimal. Other partial results for low dimensions state that biharmonic hypersurfaces with at most three distinct principal curvatures in $\R^4$ or in $\R^5$ are necessarily minimal (see \cite{Defever, HasVla95}) and very recently Fu  (see \cite{Fu2014}) extended Dimitric's result proving that any biharmonic hypersurface in $\R^n$ with at most three distinct principal curvatures is minimal. We can now state our main result:

\begin{theorem}\label{Main-theorem} Let $G$ be a cohomogeneity two group of isometries acting on $\R^n$ ($n\geq3$). Then any $G-$invariant biharmonic hypersurface in $R^n$ is minimal.
\end{theorem}

\begin{remark}The family of $G$-invariant hypersurfaces in Theorem~\ref{Main-theorem} is ample and geometrically significant (see Table~\ref{table} below). The case $d=1$ (i.e., the case of the classical rotational hypersurfaces) is a special instance of the above cited result of \cite{dimitric}. The next case, i.e., $d=2$ and $G=SO(p)\times SO(q)$, with $p+q=n$, was proved in \cite{MOR-AMPA}. For these reasons, the interest of Theorem~\ref{Main-theorem} lies on the three remaining types of $G$-actions, i.e., $d=3,\,4$ and $6$, for which the number of distinct principal curvatures is $4,\,5$ and $7$ respectively.
\end{remark}

The paper is organized as follows: in the next section we recall several basic facts from equivariant differential geometry. We believe that this short outline could be useful also for other applications in similar contexts. In the final section we provide the details of the proof of Theorem~\ref{Main-theorem}.

\section{Basic equivariant differential geometry for cohomogeneity two $G$-actions on $\R^n$}\label{equiv-diff-geom-section}

The details, together with some historical references, concerning the results of this section can be found in \cite{BackDoCarmoHsiang, HL, Pedrosa}. Let $G$ be a cohomogeneity two group of isometries acting on $\R^n$. As we mentioned in the introduction, these groups are well-understood and classified since they correspond to the \emph{isotropy representations of symmetric spaces of rank two}. For the sake of completeness, we also wish to point out the deep connection between this branch of the theory of Lie groups and the geometric properties of isoparametric functions on the Euclidean sphere (see, for example, \cite{Munzner,Thor}).

The following linear functions $w_{(d,i)}$ will play a key role:
\begin{equation}\label{definizionew_{(d,i)}}
    w_{(d,i)}(x,y)= x\, \sin (i\,\pi \slash \,d)\,-\, y\,\cos (i\,\pi \slash \,d)\,\, ,
\end{equation}
where $d$ is an integer which can be equal to $1,\,2,\,3,\,4$ or $6$, and $i$ is another integer such that $0 \leq i \leq (d-1)$. The orbit space $Q=\R^n\, / \,G$ can be identified with a linear cone of angle $(\pi \slash \,d)$ in $\R^2$ described by
\begin{equation}\label{descrizionediQ}
    Q=\left \{\, (x,y) \,\in \, \R^2 \,\, : \,\, y \geq 0 \,\,{\rm and}\,\,x\, \sin (\pi \slash \,d)\,-\, y\,\cos (\pi \slash \,d)\geq 0  \,\right \} \,\,,
\end{equation}
where the possible cases for $d$ are $1,\,2,\,3,\,4$ or $6$. We also point out that $Q=\R^2\, / \,W$, where $W=N(H,G)\, / \, H$ is the Weil group which acts on $\R^2$ by reflections with respect to the lines defined by $w_{(d,i)}(x,y)=0$. For this reason, the orbit space $Q$ is also called the \emph{Weil chamber}. The orbital distance metric $g_Q$ on $Q$ (i.e., the metric which makes the projection map $\Pi \,\, : \R^n \rightarrow Q$ a Riemannian submersion) is flat:
\begin{equation}\label{metricadiQ}
    g_Q=dx^2+dy^2
\end{equation}
and any horizontal lift of a tangent vector to $Q$ meets any $G-$orbit perpendicularly.

Let $\xi=(x,y)$ be an interior point of $Q$. We denote by $V(\xi)$ the volume of the principal orbit $\Pi^{-1}(\xi)$. The function $V(\xi)$ is called the \emph{volume function} and contains most of the information required to carry out the computation of the second fundamental form $A$ associated to a $G-$invariant hypersurface. More precisely, it turns out that $V(x,y)$ is always a homogeneous polynomial which, for each fixed type $d$, can be expressed (up to a multiplicative constant) in terms of the linear functions \eqref{definizionew_{(d,i)}} in the following form:
\begin{equation}\label{genericafunzionevolume}
    V^2(x,y)= \prod_{i=0}^{(d-1)}\,  \left [\,w_{(d,i)}(x,y)\,\right ]^{2\,m_i} \,\, ,
\end{equation}
where the $m_i$'s are positive integers (the cases which can occur are listed in Table~\ref{table}, which can be derived from an analogous Table given in \cite{Hsiang,HL}).
\begin{table}[h!]
  \centering
  \scriptsize{
\begin{tabular}{|c|c|c|c|c|}
\hline
&&&&\\
  $G$ & Action & $\dim$ Euclidean Space  & $d$ & Multiplicities \\
  &&&&\\
  \hline
  &&&& \\
  $S\mathcal{O}(n-1)$&$1+\rho_{(n-1)}$&$n\geq3$&1&$m_0=(n-2),\,m_1=1$\\
   &&&& \\
  \hline
  &&&& \\
  $S\mathcal{O}(p)\times S\mathcal{O}(q)$&$\rho_p+\rho_{q}$&$n=(p+q)\geq4$&2&$m_0=(q-1),\,m_1=(p-1)$\\
  &&&& \\
  \hline
&&&&\\
  $S\mathcal{O}(3)$ &$S^2_{\rho_{3}}-1$ &  $n=5$ & 3 & $m_0=m_1=m_2=1$ \\
  &&&&\\
  \hline
  &&&&\\
  $SU(3)$ &$Ad$ &  $n=8$ & 3 & $m_0=m_1=m_2=2$ \\
  &&&&\\
  \hline
  &&&&\\
  $Sp(3)$ & $\Lambda^2 \, \nu_{3}-1$& $ n=14$ & 3 & $m_0=m_1=m_2=4$ \\
  &&&&\\
  \hline
  &&&&\\
  $F_4$ & $\begin{array}{l}
             1 \\
             \circ \hspace{-1.3mm}- \hspace{-1.3mm}\circ \hspace{-1.4mm}=\hspace{-1.4mm}  \circ\hspace{-1.3mm} -\hspace{-1.3mm} \circ
           \end{array}
  $& $ n=26 $& 3 & $m_0=m_1=m_2=8$ \\
  &&&&\\
  \hline
  &&&&\\
  $S\mathcal{O}(5)$ & $Ad$& $ n=10$ & 4 & $m_0=m_1=m_2=m_3=2$ \\
  &&&&\\
  \hline
  &&&&\\
  $S\mathcal{O}(2)\times S\mathcal{O}(m)$ &$\rho_2 \otimes \rho_m $& $ n=2m \,\geq 6$ & 4 & $m_0=m_2=(m-2), \, m_1=m_3=1$ \\
  &&&&\\
  \hline
  &&&&\\
  $S \left ( U(2)\times U(m) \right )$ &$[\mu_2 \otimes _{\C}\mu_m]_{\R}$ & $ n=4m \,\geq 8$ & 4 & $m_0=m_2=(2m-3), \, m_1=m_3=2$ \\
  &&&&\\
  \hline
  &&&&\\
  $Sp(2) \times Sp(m)$&$ \nu_2 \otimes _{\H}\nu_m^*$& $ n=8m \,\geq 16$ & 4 & $m_0=m_2=(4m-5), \, m_1=m_3=4$ \\
  &&&&\\
  \hline
  &&&& \\
 $U(5)$ &$[\Lambda^2 \mu_5]_{\R}$&$n=20$&$4$&$m_0=m_2=5, \, m_1=m_3=4$\\
 &&&& \\
  \hline
  &&&&\\
 $ U(1)\times Spin(10)$ &$[\mu_1 \otimes _{\C}\Delta_1^+]_{\R}$ & $ n=32$ & 4 & $m_0=m_2=9, \, m_1=m_3=6$ \\
  &&&&\\
  \hline
  &&&&\\
  $G_2$ & $Ad$& $ n=14$ & 6 & $m_0=m_1= \dots =m_5=2$ \\
  &&&&\\
  \hline
  &&&&\\
  $S\mathcal{O}(4)$ &$\begin{array}{l}
              1 \hspace{1.3mm}3\\
             \circ\hspace{-1.3mm} - \hspace{-1.3mm} \circ
           \end{array} $ & $ n=8$ & 6 & $m_0=m_1= \dots =m_5=1$ \\
  &&&&\\
  \hline
  \end{tabular}
  }
  $$
  \,
  $$
  $$
  \,
  $$
  \caption{Cohomogeneity two $G$-actions on $\R^n$ (see \cite{Hsiang,HL}) (\textbf{Note}: the volume function is given in \eqref{genericafunzionevolume}, the number of distint principal curvatures of $ \Sigma_{\gamma}$ is $(d+1)$ ).}\label{table}
\end{table}

Next, we observe that any cohomogeneity one $G-$invariant immersion into $\R^n$ can be described by means of what we call its \emph{profile curve} $\gamma(s)=(x(s),\,y(s))$ in the orbit space $Q$. More precisely, the $G-$invariant hypersurface corresponding to a profile curve $\gamma$ is $\Sigma_{\gamma}= \Pi^{-1}(\gamma)$. For convenience, we shall always assume that
\begin{equation}\label{sascissacurvilinea}
    \dot{x}^2+\dot{y}^2 =1 \,\, ,
\end{equation}
and also, to fix orientation, that the unit normal $\eta$ to the hypersurface $\Sigma_{\gamma}$ projects down in $Q$ to
\begin{equation}\label{unitnormal}
    d\Pi(\eta)=\nu =-\,\dot{y}\,\frac{\partial}{\partial x}+ \dot{x}\,\frac{\partial}{\partial y}\,\, .
\end{equation}

Now, suppose that $\Sigma_{\gamma}$ is a $G-$invariant hypersurface into $\R^n$ of type $d$ ($d=1,\,2,\,3,\,4$ or $6$), with associated volume function given by \eqref{genericafunzionevolume}. Then $\Sigma_{\gamma}$ possesses $(d+1)$ distinct principal curvatures given by:
\begin{eqnarray}\label{curvatureprincipali}
    k_i &=&-\, \frac{1}{2} \, \frac{d}{d \, \nu}\, \ln \left[w_{(d,i)}(x,y) \right]^2\nonumber\\
    &=& -\, \frac{1}{2}\,  \nu\left(\ln \left[w_{(d,i)}(x,y) \right]^2\right)\,\, , \quad i=0,\, \ldots,\, (d-1) \,\,,
\end{eqnarray}
each of them with multiplicity equal to $m_i$, and
\begin{equation}\label{ultimacurvaturaprincipale}
    k_d =  \ddot{y}\,\dot {x}\,-\,\ddot{x}\,\dot {y} \,\,
\end{equation}
with multiplicity equal to one. For future use we also observe that, using \eqref{definizionew_{(d,i)}} and \eqref{unitnormal}, \eqref{curvatureprincipali} gives
\begin{eqnarray}\label{curvatureprincipali-bis}
    k_i &=& \frac{\dot{y}  \sin (i\,\pi \slash \,d)+\dot{x}  \cos (i\,\pi \slash \,d)}{w_{(d,i)}(x,y)}\nonumber\\
    &=&\frac{w_{(d,i)}(\dot{y},-\dot{x})}{w_{(d,i)}(x,y)}\,\, , \quad i=0,\, \ldots,\, (d-1) \,\,.
\end{eqnarray}
In particular, it follows that the terms $f$ and $|A|^2$ in \eqref{eq: caract_bih_normal} and \eqref{eq: caract_bih_tangent} are given by
\begin{equation}\label{espressionedif}
   f=(\, \ddot{y}\,\dot {x}\,-\,\ddot{x}\,\dot {y}\,)+\sum_{i=0}^{(d-1)}\,m_i\, k_i
\end{equation}
and
\begin{equation}\label{espressionediA^2}
    |A|^2 = (\, \ddot{y}\,\dot {x}\,-\,\ddot{x}\,\dot {y}\,)^2+ \sum_{i=0}^{(d-1)}\,m_i\, (k_i)^2
\end{equation}
respectively, where the explicit expressions for the $k_i$'s are those given in \eqref{curvatureprincipali-bis}. We shall need to compute the gradient and the laplacian of $f$: to this purpose, we work by using on $\Sigma_{\gamma}$ a local system of coordinates of type $\{u_1,\ldots,u_{n-2},s\}$, where   $u=\{u_1,\ldots,u_{n-2}\}$ are local coordinates of a principal orbit. In particular, we observe that, with respect to these local coordinates, the induced metric $g$ satisfies:
\begin{equation}\label{inducedmetric1}
  \det g = \psi(u)\, V^2(x(s),y(s))\,,
\end{equation}
where $\psi$ is a positive function on the principal orbit and
 \begin{equation}\label{inducedmetric2}
  g_{(n-1),k}=\delta_{(n-1),k}\,,\quad k=1,\ldots,\,(n-1)\,\, .
\end{equation}
Now, since $f$ depends only on $s$, it follows immediately that
\begin{equation}\label{gradientef}
\grad f = \dot {f}(s) \,\frac{\partial}{\partial s} \,\,.
\end{equation}
Next, by using \eqref{inducedmetric1} and \eqref{inducedmetric2} in
$$
\Delta f = -\,\frac{1}{\sqrt{\det g}} \,\frac{\partial}{\partial v_i}\left( g^{ij}\,\sqrt{\det g}\, \frac{\partial f}{\partial v_j}\right)\,,\quad v_i=u_i, i=1,\ldots,\,(n-2),\;\; v_{(n-1)}=s\,\,,
$$
we obtain
\begin{equation}\label{gradienteelaplaciano}
 \Delta f =- \ddot{f}- \frac{1}{2}\,\left(\frac{d}{d \, s}\, \ln V^2\right ) \, \dot{f}\,\, .
\end{equation}
We can summarize this discussion in the following proposition which, taking into account that the Ricci tensor field of $\R^n$ vanishes, follows by direct substitution of \eqref{espressionedif}, \eqref{espressionediA^2}, \eqref{gradientef} and \eqref{gradienteelaplaciano} into \eqref{eq: caract_bih_normal} and \eqref{eq: caract_bih_tangent}:

 \begin{proposition}\label{equazioniridottedibiarmonicity} Let $\Sigma_{\gamma}$ be a $G-$invariant hypersurface into $\R^n$ of type $d$ ($d=1,\,2,\,3,\,4$ or $6$), with associated volume function given by \eqref{genericafunzionevolume}. Then $\Sigma_{\gamma}$ is a biharmonic hypersurface if and only if
 \begin{equation}\label{bitensionecomponentenormale}
  \ddot{f}+ \frac{1}{2}\,\left(\frac{d}{d \, s}\, \ln V^2\right ) \, \dot{f}\, - \Big[(\, \ddot{y}\,\dot {x}\,-\,\ddot{x}\,\dot {y}\,)^2+ \sum_{i=0}^{(d-1)}\,m_i\, (k_i)^2\Big]\,f =0
 \end{equation}
and
\begin{equation}\label{bitensionecomponentetangente}
    \dot{f}\, (f+2\,(\, \ddot{y}\,\dot {x}\,-\,\ddot{x}\,\dot {y}\,)) =0 \,\, .
 \end{equation}
 \end{proposition}

 \begin{remark}\label{re:linearprifile} For future use, here we consider briefly the case where the profile curve $\gamma(s)=(x(s),y(s))$ in $Q$ satisfies $y=m x$, $m\in\R$. Therefore, we can assume that $\gamma$ is parametrized by
\begin{equation}\label{parametrizzazione}
\gamma (s)= \big( \,s \, \cos \sigma,\, s \, \sin \sigma\, \big ) \,\,,
\end{equation}
where $\sigma$ is a constant in the interval $(0,\,\pi \slash d)$. In this case, we have
\begin{equation}\label{condizioneminimalita}
    f=-\,\frac{1}{s}\,\,\sum_{i=0}^{(d-1)} \, m_i \, \cot \left (\sigma -\,\frac{i\pi}{d} \right ) \,\, .
\end{equation}
In particular, by using \eqref{parametrizzazione} and \eqref{condizioneminimalita} into \eqref{bitensionecomponentetangente}, it is immediate to deduce that a curve of this type gives rise to a solution if and only if $f \equiv 0$. In other words, for this family of curves, the associated $\Sigma_{\gamma}$ is biharmonic if and only if it is minimal.
\end{remark}

 \begin{example}\label{esempio} Equations \eqref{bitensionecomponentenormale} and \eqref{bitensionecomponentetangente} express, respectively, the vanishing of the normal and of the tangential component of the bitension field. In order to help the reader, we now compute explicitly the various principal curvatures and the biharmonicity conditions in one specific instance. More precisely, suppose that $\Sigma_{\gamma}$ is a $G-$invariant immersion into $\R^n$ of type $d=4$, with $G=U(5)$, $n=20$ (see Table~\ref{table}). In this case we have:
 $$
 w_{(4,0)}(x,y)= -y \,\, , \quad w_{(4,1)}(x,y)= \frac{\sqrt 2}{2}\,(x-y)\, ,
 $$
 $$
w_{(4,2)}(x,y)= x \,\, , \quad w_{(4,3)}(x,y)= \frac{\sqrt 2}{2}\,(x+y)\,,
 $$
 with $m_0=m_2=5$, $m_1=m_3=4$ (note that $1+\sum_{i=0}^{(d-1)}m_i= \dim (\Sigma_\gamma)=(n-1)=19$). Then we have:
$$
 \begin{aligned}
  k_0& =-\, \frac{1}{2} \, \frac{d}{d \, \nu}\, \ln \left[w_{(4,0)}(x,y)\right ]^2 \\
 &= -\, \frac{1}{w_{(4,0)}(x,y)} \left (- \dot{y}\,\frac{\partial}{\partial x}w_{(4,0)}(x,y) + \dot{x}\,\frac{\partial}{\partial y}w_{(4,0)}(x,y)\right ) = -\,\frac{\dot x}{y} \,\, ,
\end{aligned}
$$
 and, similarly,
 $$
  k_1=\frac{\dot{x}+\dot{y}}{x-y},\,\, \quad k_2= \frac{\dot y}{x}, \,\, \quad k_3=\frac{-\dot{x}+\dot{y}}{x+y} \,\, .
 $$
 Moreover, up to an irrelevant multiplicative constant,
 $$
  V^2(x,y)= \prod_{i=0}^{3}\,  \left [\,w_{(4,i)}(x,y)\,\right ]^{2\,m_i}= (xy)^{10} \, (x^2-y^2)^8 \,\, .
 $$
 Therefore, taking into account  \eqref{espressionedif}, \eqref{espressionediA^2} and \eqref{gradienteelaplaciano} the biharmonicity equations \eqref{bitensionecomponentenormale} and \eqref{bitensionecomponentetangente} become
 \begin{equation}\label{bitensionecomponentenormaleesempio}
  \ddot{f}+\,\left( 5\,\frac{\dot{x}}{x}+5\,\frac{\dot{y}}{y}+4\,\frac{\dot{x}+\dot{y}}{x+y}+4\,\frac{\dot{x}-\dot{y}}{x-y} \right ) \, \dot{f}\, - |A|^2\,f =0
 \end{equation}
and
 \begin{equation}\label{bitensionecomponentetangenteesempio}
    \dot{f}\, (f+2\,(\ddot{y}\,\dot {x}\,-\,\ddot{x}\,\dot {y})) =0
 \end{equation}
 respectively, where
 \begin{equation}\label{f-esempio}
 f= (\ddot{y}\,\dot {x}\,-\,\ddot{x}\,\dot {y})-5\, \frac{\dot x}{y}+4\,\frac{(\dot{x}+\dot{y})}{(x-y)}+5\, \frac{\dot y}{x}+4\,\frac{(-\dot{x}+\dot{y})}{(x+y)}
 \end{equation}
 and
 \begin{equation}\label{A^2esempio}
 |A|^2= (\ddot{y}\,\dot {x}\,-\,\ddot{x}\,\dot {y})^2+\,5\, \left(\frac{\dot x}{y}\right)^2+4\,\left(\frac{\dot{x}+\dot{y}}{x-y}\right)^2+5\, \left(\frac{\dot y}{x}\right)^2 +4 \,\left( \frac{-\dot{x}+\dot{y}}{x+y}\right )^2 \,\,.
 \end{equation}

 \end{example}

\section{Proof of the main Theorem}\label{proofofthemaintheorem}
\begin{proof} Now we are in the right position to prove Theorem~\ref{Main-theorem}.
The method of proof uses ideas introduced in \cite{HasVla95} and also used in \cite{MOR-AMPA}, where the case $d=2$ was proved.

It is enough to show that $\Sigma_{\gamma}$ is a CMC immersion because direct inspection of \eqref{bitensionecomponentenormale} shows that, if we have a solution with $f$ equal to a constant, then necessarily $f \equiv 0$ so that the immersion is minimal.
So, let us assume that $\Sigma_{\gamma}$ is not CMC. Then there exists a real open interval $I$ where $\dot{f}(s) \neq 0$ for all $s\in I$, and  equation \eqref{bitensionecomponentetangente} is equivalent to
\begin{equation}\label{bitensionecomponentetangentebis}
   f+2\,(\ddot{y}\,\dot {x}\,-\,\ddot{x}\,\dot {y}) =0
 \end{equation}
on $I$. Now, in order to provide a unified proof which includes all the cases (i.e., $d=1,\,2,\,3,\,4,\,6)$, it is convenient to introduce the following function:
\begin{equation}\label{definizionedellafunzioneR}
R(x,\,y,\,\dot{x},\,\dot{y})=\sum_{i=0}^{(d-1)}\,m_i\, k_i \,\, ,
\end{equation}
where the $k_i$' are those given in \eqref{curvatureprincipali-bis}. In particular, we observe from \eqref{espressionedif} that
\begin{equation}\label{legametra-f-e-R}
f= R+(\ddot{y}\,\dot {x}\,-\,\ddot{x}\,\dot {y})\,\, .
\end{equation}
Next, from equation \eqref{bitensionecomponentetangentebis} we deduce that
\begin{equation}\label{legametra-f-e-Rbis}
\ddot{y}\,\dot {x}\,-\,\ddot{x}\,\dot {y}\,=\,-\,\frac{R}{3}
\end{equation}
on $I$. Therefore, multiplying by $\dot{x}$ both sides of \eqref{legametra-f-e-Rbis} and using \eqref{sascissacurvilinea}, we easily obtain
\begin{equation}\label{eq-explicity2}
\ddot{y}= -\, \frac{R}{3}\,\,\dot{x} \,\,.
\end{equation}
In the same way, multiplying by $\dot{y}$, we also have
\begin{equation}\label{eq-explicitx2}
\ddot{x}= \frac{R}{3}\,\,\dot{y} \,\,.
\end{equation}
We also observe that, using \eqref{eq-explicity2} and \eqref{eq-explicitx2}, we can rewrite the relationship between $f$ and $R$ as follows:
\begin{equation}\label{eq-explicitf}
f=\frac{2}{3} \, R\,\,.
\end{equation}
Now, using \eqref{legametra-f-e-Rbis}$-$\eqref{eq-explicitf}, which come from the tangent component \eqref{bitensionecomponentetangente}, we claim that the normal component  \eqref{bitensionecomponentenormale} can be written in the following equivalent form:
\begin{equation}\label{eq-normal-explicit}
A_0(x,y)\, \dot{x}^3 + A_1(x,y) \,\dot{x}^2\, \dot{y} + A_2(x,y)\, \dot{x}\, \dot{y}^2
+A_3(x,y)\, \dot{y}^3=0\,\,,
\end{equation}
where the $A_j(x,y)$'s ($j=0,\, \dots,\,3$) are homogeneous polynomials of degree $3(d-1)$. In order to verify this claim, one must simply carry out an explicit computation. More precisely, first, by using \eqref{legametra-f-e-Rbis} and \eqref{eq-explicitf}, we rewrite equation \eqref{bitensionecomponentenormale} in terms of $R$ as follows:
\begin{equation}\label{bitensionecomponentenormaleinterminidiR}
  \ddot{R}+ \frac{1}{2}\,\left(\frac{d}{d \, s}\, \ln V^2\right ) \, \dot{R}\, - \Big[(1 \slash 9)\,R^2+ \sum_{i=0}^{(d-1)}\,m_i\, (k_i)^2\Big]\,R =0 \,\,.
 \end{equation}
Next, we define the following homogeneous polynomial of degree $d$ in the variables $x$, $y$:
\begin{equation}\label{definizionediQd}
    Q_d =\prod_{i=0}^{(d-1)}\, w_{(d,i)}(x,y) \,\,
\end{equation}
and also observe that
\begin{equation}\label{ddsw}
   \frac{d}{ds} w_{(d,i)}(x,y)=w_{(d,i)}(\dot{x},\dot{y})\,\, .
\end{equation}
Now, we analyse the single terms which appear in \eqref{bitensionecomponentenormaleinterminidiR}.
Using \eqref{ddsw} and \eqref{genericafunzionevolume} a direct computation gives
\begin{equation}\label{terminedivolume}
    \frac{1}{2}\,\left(\frac{d}{d \, s}\, \ln V^2\right ) =\sum_{i=0}^{(d-1)} m_i \frac{w_{(d,i)}(\dot{x},\dot{y})}{w_{(d,i)}(x,y)}= \frac{T_{1,(d-1)}\dot{x}+T_{2,(d-1)}\dot{y}}{Q_d} \,\, ,
\end{equation}
where $T_{1,(d-1)}$ and $T_{2,(d-1)}$ are the homogeneous polynomials of degree $(d-1)$ in the variables $x,\,y$ given by
$$
T_{1,(d-1)}=\sum_{i=0}^{(d-1)} m_i \frac{Q_d}{w_{(d,i)}(x,y)} \sin(i\pi/d)\,\,,
$$
$$
T_{2,(d-1)}=-\sum_{i=0}^{(d-1)} m_i \frac{Q_d}{w_{(d,i)}(x,y)} \cos(i\pi/d)\,\,.
$$
Next, using \eqref{curvatureprincipali-bis} in \eqref{definizionedellafunzioneR} we obtain
\begin{equation}\label{calcolodiR}
    R= \frac{-T_{2,(d-1)}\dot{x}+T_{1,(d-1)}\dot{y}}{Q_d}\,\, .
\end{equation}
Taking the square of \eqref{calcolodiR} and using again \eqref{curvatureprincipali-bis}, we also find
\begin{equation}\label{calcolodiA^2}
 (1 \slash 9)\,R^2+ \sum_{i=0}^{(d-1)}\,m_i\, (k_i)^2= \frac{T_{3,(2d-2)}\dot{x}^2+T_{4,(2d-2)}\dot{x}\dot{y}+T_{5,(2d-2)}\dot{y}^2}{Q_d^2}\,\,,
\end{equation}
where $T_{3,(2d-2)}$, $T_{4,(2d-2)}$  and $T_{5,(2d-2)}$ are the homogeneous polynomials of degree $(2d-2)$ in the variables $x,\,y$ given by
$$
T_{3,(2d-2)}=\frac{1}{9}\,T_{2,(d-1)}^2+\sum_{i=0}^{(d-1)} m_i \frac{Q_d^2}{[w_{(d,i)}(x,y)]^2} \cos^2(i\pi/d)
$$
$$
T_{4,(2d-2)}=-\frac{2}{9}\, T_{1,(d-1)}\,T_{2,(d-1)}+2\sum_{i=0}^{(d-1)} m_i \frac{Q_d^2}{[w_{(d,i)}(x,y)]^2} \sin(i\pi/d) \cos(i\pi/d)
$$
$$
T_{5,(2d-2)}=\frac{1}{9}\,T_{1,(d-1)}^2+\sum_{i=0}^{(d-1)} m_i \frac{Q_d^2}{[w_{(d,i)}(x,y)]^2} \sin^2(i\pi/d)\,\,.
$$
Now, taking the first and the second derivatives of \eqref{calcolodiR} and using \eqref{eq-explicity2} and \eqref{eq-explicitx2}, a direct computation shows that
\begin{equation}\label{calcolodiRpunto}
   \dot{R}= \frac{T_{6,(2d-2)}\dot{x}^2+T_{7,(2d-2)}\dot{x}\dot{y}+T_{8,(2d-2)}\dot{y}^2}{Q_d^2}\,\, ,
\end{equation}
\begin{equation}\label{calcolodiRduepunti}
   \ddot{R}= \frac{T_{9,(3d-3)}\dot{x}^3+T_{10,(3d-3)}\dot{x}^2\dot{y}+T_{11,(3d-3)}\dot{x}\dot{y}^2+T_{12,(3d-3)}\dot{y}^3}{Q_d^3}\,\, ,
\end{equation}
where $T_{6,(2d-2)},T_{7,(2d-2)},T_{8,(2d-2)}$ are homogeneous polynomials of degree $(2d-2)$ in the variables $x,\,y$, while
$T_{9,(3d-3)},T_{10,(3d-3)},T_{11,(3d-3)},T_{12,(3d-3)}$ are homogeneous polynomials of degree $(3d-3)$ in the variables $x,\,y$.

Finally, using \eqref{terminedivolume}$-$\eqref{calcolodiRduepunti} into \eqref{bitensionecomponentenormaleinterminidiR}, one obtains (up to a denominator $Q_d^3$) that the claimed equation \eqref{eq-normal-explicit} holds (note that the explicit expressions for the $A_j(x,y)$'s play no active role in the sequel, so we omit further details on this point).

Next, for a fixed $s_0\in I$, we put $x_0=x(s_0)$. Since $\dot{x}^2+\dot{y}^2=1$, we can assume that $\dot{x}(s_0)\neq 0$ and we can express $y$ as a function of $x$, $y=y(x)$, with  $x\in(x_0-\varepsilon,x_0+\varepsilon)$, and write
\begin{equation}\label{eq-non-ex-1}
\dot{y}=\frac{dy}{dx}\, \dot{x}\,\,.
\end{equation}
From $\dot{x}^2+\dot{y}^2=1$ we also obtain
\begin{equation}\label{eq-non-ex-2}
\dot{x}^2=\dfrac{1}{1+\left(\dfrac{dy}{dx}\right)^2}\,\,.
\end{equation}
For future use, combining \eqref{eq-non-ex-1} and \eqref{eq-non-ex-2}, we also observe that
\begin{equation}\label{eq-non-ex-2-bis}
\dot{x}\dot{y}=\dfrac{1}{1+\left(\dfrac{dy}{dx}\right)^2}\, \dfrac{dy}{dx}\,\,.
\end{equation}
Multiplying \eqref{eq-normal-explicit} by $\dot{x}$ and substituting \eqref{eq-non-ex-2} and \eqref{eq-non-ex-2-bis} we find that, up to a multiplicative factor $1/\left(1+\left({dy}/{dx}\right)^2\right)^2$, \eqref{eq-normal-explicit} becomes equivalent to
\begin{equation}\label{eq-non-ex-5}
A_3(x,y) \left(\dfrac{dy}{dx}\right)^3+A_2(x,y)\left(\dfrac{dy}{dx}\right)^2+A_1(x,y) \left(\dfrac{dy}{dx}\right)+A_0(x,y)=0\, \, .
\end{equation}

Next, deriving \eqref{eq-non-ex-1} with respect to $s$, a straightforward computation leads us to
\begin{equation}\label{eq-non-ex-3}
\ddot{y}=\dfrac{1}{\left(1+\left(\dfrac{dy}{dx}\right)^2\right)^2}\, \dfrac{d^2y}{dx^2}\,\,.
\end{equation}
Next, we use \eqref{eq-explicity2} into \eqref{eq-non-ex-3}: a computation, which takes into account \eqref{calcolodiR}, \eqref{eq-non-ex-2} and \eqref{eq-non-ex-2-bis}, shows that the expression for the second derivative of $y$ with respect to $x$ is given by:
\begin{equation}\label{eq-non-ex-4}
\begin{aligned}
\frac{d^2y}{dx^2}=& -\frac{1}{3} \, \left(1+\left(\dfrac{dy}{dx}\right)^2\right)^2\, R\, \dot{x}\\
=& -\frac{1}{3} \, \left(1+\left(\dfrac{dy}{dx}\right)^2\right)^2\, \left(\frac{-T_{2,(d-1)}\dot{x}+T_{1,(d-1)}\dot{y}}{Q_d}\right)\, \dot{x}\\
=& \frac{1}{3}\, \left(
1+\left(\dfrac{dy}{dx}\right)^2\right)\, \left (\frac{ T_{2,(d-1)}-T_{1,(d-1)}\,(dy \slash dx) }{Q_d} \right )\,\,.
\end{aligned}
\end{equation}

Next, taking the derivative of \eqref{eq-non-ex-5} with respect to $x$, that is $d/dx$, and using \eqref{eq-non-ex-4}, we obtain, up to $1/Q_d$, the following equation:
\begin{eqnarray}\label{eq-non-ex-6}
C_5(x,y) \left(\dfrac{dy}{dx}\right)^5&+&C_4 (x,y) \left(\dfrac{dy}{dx}\right)^4 \\ \nonumber
+\,\, C_3(x,y)\left(\dfrac{dy}{dx}\right)^3&+&C_2(x,y)
\left(\dfrac{dy}{dx}\right)^2+\,C_1(x,y)\left(\dfrac{dy}{dx}\right)+\, C_0(x,y)=0\,\,,
\end{eqnarray}
where the $C_j(x,y)$'s ($j=0,\, \dots,\,5$) are homogeneous polynomials of degree $4(d-1)$ which are related to the $A_j(x,y)$'s as follows:
\begin{equation}\label{definizionedeipolinomiC}
    \left \{\begin{array}{l}
           C_0=Q_d \cdot {\partial A_{0}}/{\partial x}+(1/3)\,A_1 \cdot T_{2,(d-1)}  \\
           C_1=Q_d \cdot ({\partial A_{1}}/{\partial x}+{\partial A_{0}}/{\partial y})+(2/3)\,A_2 \cdot T_{2,(d-1)}-(1/3)\, A_1 \cdot T_{1,(d-1)} \\
           C_2=Q_d \cdot ({\partial A_{2}}/{\partial x}+{\partial A_{1}}/{\partial y})+A_3 \cdot T_{2,(d-1)}-(2/3)\,A_2 \cdot T_{1,(d-1)}+(1/3)\,A_1 \cdot T_{2,(d-1)} \\
           C_3=Q_d \cdot ({\partial A_{3}}/{\partial x}+{\partial A_{2}}/{\partial y})-A_3 \cdot T_{1,(d-1)}+(2/3)\,A_2 \cdot T_{2,(d-1)}-(1/3)\, A_1 \cdot T_{1,(d-1)}\\
           C_4=Q_d \cdot {\partial A_{3}}/{\partial y}+A_3 \cdot T_{2,(d-1)}-(2/3)\,A_2\cdot T_{1,(d-1)} \\
           C_5=-\,A_3 \cdot T_{1,(d-1)} \,\,,
           \end{array}
    \right .
\end{equation}
where in this computation we have used 
$$
\frac{d A_i}{dx}=\frac{\partial A_{i}}{\partial x}+\frac{\partial A_{i}}{\partial y } \dfrac{dy}{dx}\,,\quad i=0,\ldots,3\,.
$$

Now we are in the right position to end the proof: for any arbitrarily fixed $x_1\in(x_0-\varepsilon,x_0+\varepsilon)$, setting $y_1=y(x_1)$, \eqref{eq-non-ex-5} and \eqref{eq-non-ex-6}
can be thought of as two polynomial equations in $dy/dx$, with coefficients given, respectively,  by $A_j(x_1,y_1),\, j=0,\ldots,3$ and $C_j(x_1,y_1),\, j=0,\ldots,5$,  which have the common solution $(dy/dx)(x_1)$.
Using standard arguments of algebraic geometry (\cite{HasVla95}), this implies that the resultant of the two polynomials is zero for any $x_1\in(x_0-\varepsilon,x_0+\varepsilon)$.  Now, since the coefficients $A_j(x,y)$ and $C_j(x,y)$ are {\it homogeneous
polynomials} (of degree $3(d-1)$ and $4(d-1)$ respectively), it turns out that the resultant is itself a {\it homogeneous polynomial} of degree $3(d-1)\cdot 5+ 4(d-1)\cdot 3=27(d-1)$. Then, we can divide the resultant by $x^{27(d-1)}$ and putting
$z=y/x$ we obtain a polinomial equation in $z$ with constant coefficients. Since  
$z$ is continuous it must be constant, that is $y=mx$, $m\in\R$.  But any such solution is biharmonic if and only if it is minimal (see Remark~\ref{re:linearprifile}): a contradiction with the hypothesis that $\Sigma_{\gamma}$ is not CMC. 
\end{proof}
\begin{remark} We think that it is important to stress the fact the method of proof works because the $G$-invariance is sufficient to guarantee that the biharmonicity conditions \eqref{bitensionecomponentenormale} and \eqref{bitensionecomponentetangente}, when expressed in the form \eqref{eq-non-ex-5} and \eqref{eq-non-ex-6}, have coefficients given by \emph{homogeneous} polynomials. In the direction of proving the Chen conjecture, we also point out that the incompatibility between \eqref{bitensionecomponentenormale} and \eqref{bitensionecomponentetangente} appear to be of a local nature.
\end{remark}

\begin{example} For the sake of completeness, we report here the explicit expressions of the relevant homogeneous polynomials which appear in the Example~\ref{esempio}. In this case we have:

$$
A_1(x,y)=A_2(y,x)=8 \left(-25 x^2 y^7+540 x^4 y^5-118 x^6 y^3+65 x^8
   y-10 y^9\right)\,,
$$

$$
A_0(x,y)=A_3(y,x)=775 x^7 y^2-363 x^5 y^4+1237 x^3 y^6+130 x y^8-275
   x^9\,.
$$

\end{example}

\end{document}